\documentclass[a4paper,11pt]{amsart}

\usepackage[T1]{fontenc}
\usepackage{lmodern}
\usepackage{amsmath,amssymb,amsthm}
\usepackage{enumitem,needspace}
\usepackage{microtype}
\usepackage[margin=28mm,headheight=14pt]{geometry}
\usepackage[hidelinks]{hyperref}

\makeatletter
\renewcommand{\@cite}[2]{{\@citestyle[\citeform{#1}\if@tempswa, #2\fi]}}
\renewcommand{\@biblabel}[1]{[#1]}
\makeatother

\newlist{claimparts}{enumerate}{1}
\setlist[claimparts]{label=(\arabic*),ref=\arabic*,leftmargin=2.5em,
  labelsep=.6em,itemsep=3pt,topsep=4pt,parsep=0pt}
\numberwithin{equation}{section}
\theoremstyle{plain}
\newtheorem{theorem}{Theorem}[section]
\newtheorem{lemma}[theorem]{Lemma}
\newtheorem{proposition}[theorem]{Proposition}

\theoremstyle{remark} 
\newtheorem{remark}[theorem]{Remark} 

\DeclareMathOperator{\Div}{\mathrm{div}}
\DeclareMathOperator{\tr}{\mathrm{Tr}}

\allowdisplaybreaks[2]
\title[Curvature pinching in spheres]{Curvature pinching for minimal Tangent-Yang-Mills hypersurfaces in spheres}

\author[J.Q. Ge]{Jianquan Ge}
\address{School of Mathematical Sciences, Laboratory of Mathematics and Complex Systems, Beijing Normal University, Beijing 100875, P. R. China}
\email{jqge@bnu.edu.cn}

\author[L.X. Xiao]{Lixin Xiao\textsuperscript{*}}
\address{Chern Institute of Mathematics and LPMC, Nankai University, Tianjin 300071, P. R. China}
\email{lixinxiao@nankai.edu.cn}

\author[W.J. Zhang]{Wenjin Zhang}
\address{Chern Institute of Mathematics and LPMC, Nankai University, Tianjin 300071, P. R. China}
\email{wjzhang@nankai.edu.cn}

\date{}
\thanks{\textsuperscript{*}Corresponding author.}
\thanks{J. Q. Ge is partially supported by NSFC (No. 12571049) and the Fundamental Research Funds for the Central Universities.}

\begin{document}

\begin{abstract}
We prove a rigidity theorem for closed connected minimal Tangent-Yang-Mills hypersurfaces in the unit sphere. In every dimension $n\geqslant 3$, such a hypersurface is totally geodesic if the squared norm of its Riemannian curvature
tensor is at most $2n(n-1)$.
\end{abstract}
\maketitle

\section{Introduction}
A classical problem for minimal hypersurfaces in spheres is to determine which curvature pinching conditions force rigidity. Let $M^n\to\mathbb S^{n+1}(1)$ be a closed connected minimal hypersurface, where $n\geqslant 3$. We write $h$ for its second fundamental form and $A$ for its shape operator, and set
$S=|h|^2=\tr(A^2)$. The classical pinching theorem of Simons, Chern-do Carmo-Kobayashi, and Lawson \cite{ChernDoCarmoKobayashi,Lawson,Simons} states that $S\leqslant  n$ forces $M$ to be totally geodesic or a minimal Clifford hypersurface. Thus the first positive constant value of $S$ is $n$, attained precisely by the minimal Clifford hypersurfaces.

This first gap leads to the question of which  constant values of $S$ can occur beyond $n$, a problem that forms part  of Chern's conjecture. The second-gap problem asks for an interval above $n$ containing no further constant values of $S$; the initial estimates in this direction were obtained
by Peng and Terng \cite{PengTerng}. Recently, Cheng, Li, Wei, and Yin \cite{ChengLiWeiYin} proved that a complete minimal hypersurface with constant scalar curvature and $n\geqslant4$ satisfies $S>\frac32n$ whenever $S>n$. In the pointwise pinching problem, Ding and Xin \cite{DingXin} proved that there exists $\delta(n)>0$ such that $n\leqslant  S\leqslant  n+\delta(n)$ forces $M$ to be a minimal Clifford hypersurface,
without assuming $S$ to be constant. For $n\geqslant 6$, they obtained $\delta(n)=\frac{n}{23}$; subsequent refinements appear in \cite{LeiXuXu, XuXu}.

Further rigidity results are available under additional assumptions on intrinsic curvature or on the extrinsic curvature functions $f_r=\tr(A^r)$, $r\geqslant 1$. Cheng and Ishikawa \cite{ChengIshikawa} proved that a closed minimal
hypersurface with Ricci curvature bounded below by $\frac{n}{2}$ and $n\leqslant  S\leqslant  n+\frac{14(n+4)}{9n+30}$ is a Clifford hypersurface. For closed oriented minimal hypersurfaces with constant $S$ and $f_3$, Ge, Tan, Yan, and Zhang \cite{GeTanYanZhang} proved that $S>n$ implies $S\geqslant2n$, with equality precisely for the minimal Cartan hypersurfaces with three principal curvatures.
These results distinguish pointwise pinching in the general minimal class from larger gaps under additional geometric conditions.

In this paper, we study rigidity under a pointwise upper bound
on the squared norm of the Riemannian curvature tensor.
This squared norm also enters the Yang-Mills functional
of the Levi-Civita connection.
Let $R_{ijkl}$ be its components and set
\[
R=\frac{1}{2}\sum_{i,j,k,l}R_{ijkl}^2.
\]
The factor $\frac{1}{2}$ reflects the standard normalization for the norm of curvature viewed as a two-form with values in skew-adjoint endomorphisms.
For a minimal hypersurface, the Gauss equation gives
\[
R=n(n-1)+S^2-2S-f_4.
\]
Thus a pointwise upper bound on $R$ imposes a constraint involving both $S$ and $f_4$. A totally geodesic hypersurface has $S=0$ and $R=n(n-1)$, but, unlike $S=0$, the bound $R\leqslant  n(n-1)$ does not characterize total geodesicity. Indeed, the minimal Clifford hypersurface
\[
\mathbb S^1\!\left(\sqrt{\frac{1}{n}}\right)
\times\mathbb S^{n-1}\!\left(\sqrt{\frac{n-1}{n}}\right)
\subset\mathbb S^{n+1}(1)
\]
has $R=\frac{n^2(n-2)}{n-1}<n(n-1)$.

This example indicates that an additional condition is necessary for rigidity at the totally geodesic value of $R$. The variational problem for the same curvature density provides such a condition. Ge, Xiao, and Zhang \cite{GeXiaoZhang} call an immersion Tangent-Yang-Mills if it is a
critical point of
\[
\int_M R\,\mathrm{dvol}
\]
under arbitrary normal variations of the immersion, with both the induced metric and its Levi-Civita connection varying accordingly. This condition relates the extrinsic curvature functions and their derivatives. For a minimal hypersurface in $\mathbb S^{n+1}(1)$, the Tangent-Yang-Mills equation takes the form
\begin{equation}\label{TangentYangMillsEquation}
\frac{1}{2}\Div(A\nabla S)-\frac{1}{3}\Delta f_3
+(S-2)f_3-f_5=0.
\end{equation}

We prove that the pointwise bound of $R$ by $n(n-1)$
forces total geodesicity.
\begin{theorem}\label{MainTheorem}
Let $M^n\to\mathbb S^{n+1}(1)$ be a closed connected minimal
Tangent-Yang-Mills hypersurface, where $n\geqslant 3$. If
\[
R\leqslant  n(n-1)
\]
pointwise on $M$, then $M$ is totally geodesic.
\end{theorem}

\begin{remark}
For $n\geqslant 4$, the estimates in Lemma~\ref{CompressionMomentInequalities}
give $f_4\leqslant \frac{n^2-3n+3}{n(n-1)}S^2$. Hence   our bound $R\leqslant  n(n-1)$ implies
\[
S\leqslant \frac{2n(n-1)}{2n-3}=n+\frac{n}{2n-3}.
\]
However, the converse need not hold. In $\mathbb S^5(1)$, the  minimal Clifford hypersurfaces
  \[
\mathbb S^1\!\left(\frac{1}{2}\right)\times
\mathbb S^3\!\left(\frac{\sqrt3}{2}\right),\qquad
\mathbb S^2\!\left(\frac{1}{\sqrt2}\right)\times
\mathbb S^2\!\left(\frac{1}{\sqrt2}\right)
\]
both have $S=4$, whereas their values of $R$ are
$\frac{32}{3}$ and 16, respectively.
In dimension three, $f_4=\frac{S^2}{2}$, so the curvature hypothesis is equivalent to $S\leqslant 4$.
The resulting upper bound on $S$ exceeds the classical threshold $S=n$ in every dimension $n\geqslant 3$. 
\end{remark}

Section~\ref{Preliminaries} derives an integral identity relating the extrinsic curvature functions to their gradients. For $n\geqslant4$, Section~\ref{HigherDimensions} estimates its two sides and compares them with a weighted Simons identity to show that $S$ is constant; the resulting
rigidity then follows from the curvature bound and the Tangent-Yang-Mills condition. Section~\ref{DimensionThree} proves the result for $n=3$ by introducing an additional weighted integral identity and decomposing the resulting integrand into two nonnegative terms whose simultaneous
vanishing forces total geodesicity.

\section{Preliminaries}\label{Preliminaries}

Let $M^n\to\mathbb S^{n+1}(1)$ be a minimal hypersurface with $n\geqslant 3$. We write $A$ for its shape operator and $h$ for its second fundamental form, and use $\nabla$ for the induced covariant derivative. At a fixed point, choose a local orthonormal frame $e_1,\ldots,e_n$ that diagonalizes $A$,
and denote its principal curvatures by $\lambda_1,\ldots,\lambda_n$. Thus
\[
Ae_i=\lambda_i e_i,\qquad h_{ij}=h(e_i,e_j)=\lambda_i\delta_{ij},\qquad
f_r=\tr(A^r)=\sum_i\lambda_i^r.
\]
Here $r$ is a positive integer, and $f_2=S=|h|^2$.
All norms and inner products are induced by the metric, and $\mathrm{dvol}$ denotes its volume element.

We set $\nabla_k=\nabla_{e_k}$ and use subscripts for covariant derivatives.
For scalar functions, a subscript $k$ denotes differentiation in the direction $e_k$; in particular, $S_k=\nabla_kS$. Thus
\[
h_{ijk}=(\nabla_{e_k}h)(e_i,e_j),\qquad (f_r)_k=\nabla_kf_r.
\]
Minimality and the Codazzi equations give
\[
\sum_i\lambda_i=0,\qquad
h_{ijk}=h_{jik}=h_{ikj},\qquad
\sum_i h_{iik}=0.
\]
Differentiating the trace of $A^r$ at this point gives
\begin{equation}\label{CurvatureMomentDerivatives}
(f_r)_k=r\sum_i\lambda_i^{r-1}h_{iik}\qquad(r\geqslant 2).
\end{equation}

The Gauss equation reads
\[
R_{ijkl}=\delta_{ik}\delta_{jl}-\delta_{il}\delta_{jk}
+h_{ik}h_{jl}-h_{il}h_{jk}.
\]
It follows that
\begin{align*}
R=\frac{1}{2}\sum_{i,j,k,l}R_{ijkl}^2=2\sum_{i<j}(1+\lambda_i\lambda_j)^2=n(n-1)+S^2-2S-f_4.
\end{align*}
Thus the curvature hypothesis is equivalent to the following bound for $f_4$:
\begin{equation}\label{CurvatureMomentBound}
R\leqslant  n(n-1)\quad\Longleftrightarrow\quad S^2-2S-f_4\leqslant 0
\quad\Longleftrightarrow\quad f_4\geqslant  S^2-2S.
\end{equation}
Suppose now that $M$ is closed and satisfies the Tangent-Yang-Mills equation
\eqref{TangentYangMillsEquation}. To compare its algebraic and derivative
terms, multiply by $f_3$ and integrate:
\[
0=\int_M\left(\frac{1}{2}f_3\Div(A\nabla S)
-\frac{1}{3}f_3\Delta f_3+(S-2)f_3^2-f_3f_5\right)\,\mathrm{dvol}.
\]
Since $M$ is closed, integration by parts gives
\begin{align*}
0&=\int_M\Div(f_3A\nabla S)\,\mathrm{dvol}\\
 &=\int_M f_3\Div(A\nabla S)\,\mathrm{dvol}
   +\int_M\langle A\nabla S,\nabla f_3\rangle\,\mathrm{dvol},\\
0&=\int_M\Div(f_3\nabla f_3)\,\mathrm{dvol}\\
 &=\int_M f_3\Delta f_3\,\mathrm{dvol}
   +\int_M|\nabla f_3|^2\,\mathrm{dvol}.
\end{align*} 
Substituting these identities separates the zeroth-order terms from the derivative terms:
\[
0=\int_M\left((S-2)f_3^2-f_3f_5\right)\,\mathrm{dvol}
+\int_M\left(\frac{1}{3}|\nabla f_3|^2
-\frac{1}{2}\langle A\nabla S,\nabla f_3\rangle\right)\,\mathrm{dvol}.
\]
Define
\begin{equation}\label{AlgebraicIntegrandDefinition}
\mathcal P=f_3f_5-(S-2)f_3^2,
\end{equation}
\begin{equation}\label{DerivativeIntegrandDefinition}
\mathcal Q=\frac{1}{3}|\nabla f_3|^2
 -\frac{1}{2}\langle A\nabla S,\nabla f_3\rangle.
\end{equation}
With this notation, the preceding integral identity becomes
  \begin{equation}\label{CubicMomentIntegralSplit}
\begin{aligned}
\int_M\mathcal P\,\mathrm{dvol}
 &=\int_M\left(f_3f_5-(S-2)f_3^2\right)\,\mathrm{dvol}\\
 &=\int_M\left(\frac{1}{3}|\nabla f_3|^2
  -\frac{1}{2}\langle A\nabla S,\nabla f_3\rangle\right)\,\mathrm{dvol}
 =\int_M\mathcal Q\,\mathrm{dvol}.
\end{aligned}
\end{equation} 

\section{The case \texorpdfstring{$n\geqslant 4$}{n >= 4}}\label{HigherDimensions}

For a closed minimal Tangent-Yang-Mills hypersurface of dimension $n\geqslant4$, we begin with the expression
$\mathcal Q$ in \eqref{DerivativeIntegrandDefinition}. Completing the square gives
\[
\mathcal Q=3\left|\frac13\nabla f_3-\frac14A\nabla S\right|^2
-\frac3{16}|A\nabla S|^2=3|T|^2-\frac3{16}|A\nabla S|^2,
\]
where
\begin{equation}\label{CompletedGradientDefinition}
T=\frac13\nabla f_3-\frac14A\nabla S.
\end{equation}
In a principal frame, \eqref{CurvatureMomentDerivatives} gives
\[
T_k=\sum_i\left(\lambda_i^2-\frac{\lambda_k\lambda_i}{2}\right)h_{iik}.
\]
The coefficients satisfy
\[
\left(\sum_i\left(\lambda_i^2-\frac{\lambda_k}{2}\lambda_i\right)^2\right)^{\frac{1}{2}}=\left|A^2-\frac{\lambda_k}{2}A\right|
\leqslant|A^2|+\frac{|\lambda_k|}{2}|A|
\leqslant\frac32S.\]
Consequently, Cauchy-Schwarz yields
\begin{equation}\label{CoarseCompletedGradientEstimate}
|T|^2\leqslant\frac94S^2\sum_{i,k}h_{iik}^2
\leqslant\frac94S^2|\nabla h|^2,
\qquad
\mathcal Q\leqslant\frac{27}{4}S^2|\nabla h|^2.
\end{equation}
This preliminary estimate suggests comparing the Tangent-Yang-Mills integral identity~\eqref{CubicMomentIntegralSplit} with a weighted
form of the following identity of Simons \cite{Simons}:
\begin{equation}\label{SimonsIdentity}
\frac12\Delta S=|\nabla h|^2+S(n-S).
\end{equation}
Multiplying by $S^2$ and integrating by parts gives
\[
\begin{aligned}
\int_M\left(S^2|\nabla h|^2+S^3(n-S)\right)\,\mathrm{dvol}
 =\frac{1}{2}\int_M S^2\Delta S\,\mathrm{dvol}=-\int_M S|\nabla S|^2\,\mathrm{dvol}.
\end{aligned}
\]
Hence
\begin{equation}\label{WeightedSimonsIntegral}
\int_M S^3(S-n)\,\mathrm{dvol}
=\int_M\bigl(S^2|\nabla h|^2+S|\nabla S|^2\bigr)\,\mathrm{dvol}.
\end{equation}
We compare \eqref{CubicMomentIntegralSplit} with
\eqref{WeightedSimonsIntegral} by estimating the algebraic terms on the left and the gradient terms on the right. Propositions~\ref{CurvatureMomentLowerEstimate}
and~\ref{GradientUpperEstimate} establish the pointwise bounds
\begin{equation}\label{AlgebraicIntegrandTargetBound}
\mathcal P\geqslant \frac{2(n-2)^2}{n^2}S^3(S-n),
\end{equation}
\begin{equation}\label{DerivativeIntegrandTargetBound}
\mathcal Q\leqslant \frac{3(5n-8)}{14(n+2)}
 \bigl(S^2|\nabla h|^2+S|\nabla S|^2\bigr).
\end{equation}
Once these pointwise estimates are established, we obtain
\begin{equation}\label{HigherDimensionalIntegralComparison}
\begin{aligned}
\frac{2(n-2)^2}{n^2}
 \int_M\bigl(S^2|\nabla h|^2+S|\nabla S|^2\bigr)\,\mathrm{dvol}
&=\frac{2(n-2)^2}{n^2}\int_M S^3(S-n)\,\mathrm{dvol}\\
&\leqslant \int_M\mathcal P\,\mathrm{dvol}=\int_M\mathcal Q\,\mathrm{dvol}\\
&\leqslant \frac{3(5n-8)}{14(n+2)}
 \int_M\bigl(S^2|\nabla h|^2+S|\nabla S|^2\bigr)\,\mathrm{dvol}.
\end{aligned}
\end{equation}
For $n\geqslant4$, the coefficients satisfy
\begin{equation}\label{HigherDimensionalCoefficientGap}
\begin{aligned}
\frac{2(n-2)^2}{n^2}-\frac{3(5n-8)}{14(n+2)}
 =\frac{1}{14}
 +\frac{(n-4)\left(6n^2+7(n-4)\right)}{7n^2(n+2)}>0.
\end{aligned}
\end{equation}
Thus the two estimates force the integrals in
\eqref{WeightedSimonsIntegral} to vanish. It follows that $S$ is constant
and equals either 0 or $n$. In Theorem~\ref{HigherDimensionalRigidity},
we complete the proof by ruling out $S\equiv n$.

To prove \eqref{AlgebraicIntegrandTargetBound} and
\eqref{DerivativeIntegrandTargetBound}, we first derive algebraic relations among the curvature functions $f_r$ by computing the traces of powers of a compression of $A$. Choose a principal orthonormal frame $\{e_1,\ldots,e_n\}$
and set
\begin{equation}\label{PrincipalFrameVectors}
\eta_1=\sum_i e_i,\qquad \eta_2=\sum_i\lambda_i e_i,
\qquad V=\operatorname{span}\{\eta_1,\eta_2\}.
\end{equation}
By definition,
\[
|\eta_1|^2=n,\qquad |\eta_2|^2=S,\qquad A\eta_1=\eta_2.
\]
Minimality further gives
\[
\langle\eta_1,\eta_2\rangle=\sum_i\lambda_i=0.
\]
Thus, at a point where $S>0$, the vectors $\eta_1$ and $\eta_2$ are nonzero
and orthogonal.
For every integer $r\geqslant 1$, we have
\[
A^r\eta_1=\sum_i\lambda_i^r e_i,
\]
and hence
\[
\langle A^r\eta_1,\eta_1\rangle=f_r,\qquad
\langle A^r\eta_1,\eta_2\rangle=f_{r+1},\qquad
\langle A^r\eta_2,\eta_2\rangle=f_{r+2}.
\]

Let $\pi$ be the orthogonal projection onto $V^\perp$. Since $\eta_1$ and $\eta_2$ are orthogonal,
\[
\pi(v)=v-\frac{\langle v,\eta_1\rangle}{n}\eta_1
          -\frac{\langle v,\eta_2\rangle}{S}\eta_2.
\]
In particular, $\pi^2=\pi$ and $\pi$ is self-adjoint. Define the compression of the shape operator $A$ to $V^\perp$ by
\[
A^{\perp}=(\pi\circ A)|_{V^\perp}:V^\perp\longrightarrow V^\perp.
\]
  Since  $\pi$  is an orthogonal projection  and $A$   is self-adjoint,
$\langle A^{\perp}v,w\rangle=\langle Av,w\rangle=\langle v,Aw\rangle=\langle v,A^{\perp}w\rangle$
for $v,w\in V^\perp$, so $A^{\perp}$ is self-adjoint.
We next compute the traces of its first three powers.

\begin{lemma}\label{CompressionTraceIdentities}
For   $n\geqslant 3$, at every point where $S>0$, the following identities hold:
\begin{claimparts}
\item \(\displaystyle \tr A^{\perp}=-\frac{f_3}{S}\);
\item \(\displaystyle \tr (A^{\perp})^2=S-\frac{2f_4}{S}+\frac{f_3^2}{S^2}\);
\item \(\displaystyle \tr (A^{\perp})^3=f_3-\frac{3f_5}{S}
 +\frac{3f_3f_4}{S^2}-\frac{f_3^3}{S^3}\).
\end{claimparts}
\end{lemma}
\begin{proof}
Choose an orthonormal basis $\widetilde\eta_1,\ldots,\widetilde\eta_n$ with
\[
\widetilde\eta_1=\frac{\eta_1}{\sqrt n},\qquad
\widetilde\eta_2=\frac{\eta_2}{\sqrt S}.
\]
In this basis,
\[
A=\begin{pmatrix}A_1&A_2\\A_2^{\mathrm T}&A_4\end{pmatrix},
\qquad\pi=\begin{pmatrix}0&0\\0&I_{n-2}\end{pmatrix},
\qquad A^{\perp}=A_4,
\]
where
\[
A_1\in M_{2\times2}(\mathbb R),\qquad
A_2\in M_{2\times(n-2)}(\mathbb R),\qquad
A_4\in M_{(n-2)\times(n-2)}(\mathbb R).
\]
Here $M_{r\times s}(\mathbb R)$ is the space of real $r\times s$ matrices,
and $I_{n-2}$ is the identity on $V^\perp$.
Computing the entries of the first block, we obtain
\[
A_1=\begin{pmatrix}
\langle A\widetilde\eta_1,\widetilde\eta_1\rangle&
\langle A\widetilde\eta_2,\widetilde\eta_1\rangle\\[2pt]
\langle A\widetilde\eta_1,\widetilde\eta_2\rangle&
\langle A\widetilde\eta_2,\widetilde\eta_2\rangle
\end{pmatrix}
=\begin{pmatrix}
0&\sqrt{\frac{S}n}\\[2pt]
\sqrt{\frac{S}n}&\frac{f_3}{S}
\end{pmatrix}.
\]

The operator $\pi\circ A\circ\pi$ agrees with $A^{\perp}$ on $V^\perp$ and vanishes on $V$. Since $\pi^2=\pi$, cyclicity of the trace gives
\[
\tr (A^{\perp})^r=\tr\Bigl((\pi\circ A\circ\pi)^r\Bigr)
       =\tr\Bigl((A\pi)^r\Bigr)
       =\tr\Bigl(\bigl(A-A(I-\pi)\bigr)^r\Bigr),\qquad r\geqslant1.
\]
For the first trace, this yields
\[
\tr A^{\perp}=\tr A-\tr\Bigl(A(I-\pi)\Bigr)
      =-\tr A_1=-\frac{f_3}{S}.
\]
For the second trace, we similarly obtain
\begin{align*}
\tr (A^{\perp})^2=\tr A^2-2\tr\Bigl(A^2(I-\pi)\Bigr)
  +\tr\Bigl(\bigl(A(I-\pi)\bigr)^2\Bigr),
\end{align*}
where
\begin{align*}
\tr\Bigl(A^2(I-\pi)\Bigr)
&=\langle A^2\widetilde\eta_1,\widetilde\eta_1\rangle
 +\langle A^2\widetilde\eta_2,\widetilde\eta_2\rangle
 =\frac{S}n+\frac{f_4}{S},\\
\tr\Bigl((A(I-\pi))^2\Bigr)
&=\tr A_1^2=\frac{2S}{n}+\frac{f_3^2}{S^2}.
\end{align*}
Consequently,
\[
\tr (A^{\perp})^2
=S-2\left(\frac{S}n+\frac{f_4}{S}\right)
 +\frac{2S}{n}+\frac{f_3^2}{S^2}
=S-\frac{2f_4}{S}+\frac{f_3^2}{S^2}.
\]
Similarly, expanding the cube gives
\begin{align*}
\tr (A^{\perp})^3
&=f_3-3\left(\frac{f_3}{n}+\frac{f_5}{S}\right)
 +3\left(\frac{2f_3}{n}+\frac{f_3f_4}{S^2}\right)
 -\left(\frac{3f_3}{n}+\frac{f_3^3}{S^3}\right)\\
&=f_3-\frac{3f_5}{S}+\frac{3f_3f_4}{S^2}-\frac{f_3^3}{S^3}.
\end{align*} 
\end{proof}

We now define the traceless part of $A^{\perp}$ by
\[
\mathring{A}^{\perp}=A^{\perp}-\frac{\tr A^{\perp}}{n-2}I_{n-2}
           =A^{\perp}+\frac{f_3}{(n-2)S}I_{n-2}.
\]
It is self-adjoint and traceless. In the next lemma, we compute
$\tr(\mathring{A}^{\perp})^2$ to obtain an upper bound for $f_4$ and $\tr(\mathring{A}^{\perp})^3$ to solve for $f_5$.
Proofs of inequalities~(1) and~(2) are also given in
\cite{SharmaBhandari} and \cite{Okumura}, respectively.

\begin{lemma}\label{CompressionMomentInequalities}
For $n\geqslant 3$, at every point where $S>0$, the following assertions hold:
\begin{claimparts}
\item \(\displaystyle f_4\leqslant \frac{S^2}{2}
 +\frac{n-3}{2(n-2)}\frac{f_3^2}{S}\);
\item \(\displaystyle f_3^2\leqslant \frac{(n-2)^2}{n(n-1)}S^3\);
\item \(\displaystyle \frac{f_3^2}{S^2}+\tr(\mathring{A}^{\perp})^2
 \leqslant \frac{n-2}{n}S\);
\item
\(
\tr(\mathring{A}^{\perp})^3
 =\dfrac{n+1}{n-2}f_3
 +\dfrac{3(n-4)}{n-2}\dfrac{f_3f_4}{S^2}-\dfrac{(n-3)(n-4)}{(n-2)^2}\dfrac{f_3^3}{S^3}
 -\dfrac{3f_5}{S};
\)
\item \(\displaystyle |f_3\tr(\mathring{A}^{\perp})^3|
 \leqslant \frac{n-2}{2n}S^2\tr(\mathring{A}^{\perp})^2\).
\end{claimparts}
\end{lemma}

\begin{proof}
Expanding the square of $\mathring{A}^{\perp}$ and applying
Lemma~\ref{CompressionTraceIdentities}, we obtain
\begin{align*}
\tr(\mathring{A}^{\perp})^2
&=\tr (A^{\perp})^2+\frac{2f_3}{(n-2)S}\tr A^{\perp}
 +\frac{f_3^2}{(n-2)S^2}\\
&=S-\frac{2f_4}{S}+\frac{n-3}{n-2}\frac{f_3^2}{S^2}.
\end{align*}
Since $\mathring{A}^{\perp}$ is self-adjoint, $\tr(\mathring{A}^{\perp})^2\geqslant0$.
Consequently,
\begin{equation}\label{FourthMomentCompressionBound}
\begin{aligned}
f_4&=\frac{1}{2}\left(S^2+\frac{n-3}{n-2}\frac{f_3^2}{S}
 -S\tr(\mathring{A}^{\perp})^2\right)\leqslant\frac{S^2}{2}+\frac{n-3}{2(n-2)}\frac{f_3^2}{S},
\end{aligned}
\end{equation}
which proves (1). Similarly, expanding the cube gives (4):
  \begin{equation}\label{CenteredCompressionCube}
\begin{aligned}
\tr(\mathring{A}^{\perp})^3
&=\tr (A^{\perp})^3+\frac{3f_3}{(n-2)S}\tr (A^{\perp})^2
 +\frac{3f_3^2}{(n-2)^2S^2}\tr A^{\perp}
 +\frac{f_3^3}{(n-2)^2S^3}\\
&=\frac{n+1}{n-2}f_3
 +\frac{3(n-4)}{n-2}\frac{f_3f_4}{S^2}-\frac{(n-3)(n-4)}{(n-2)^2}\frac{f_3^3}{S^3}
 -\frac{3f_5}{S}.
\end{aligned}
\end{equation} 

To prove (2) and (3), we decompose $A\eta_2$ orthogonally:
\[
A\eta_2=\frac{S}{n}\eta_1+\frac{f_3}{S}\eta_2+\pi(A\eta_2),
\]
since $\langle A\eta_2,\eta_1\rangle=S$ and
$\langle A\eta_2,\eta_2\rangle=f_3$.
Taking squared norms gives
\[
f_4=|A\eta_2|^2
=\frac{S^2}{n}+\frac{f_3^2}{S}+|\pi(A\eta_2)|^2.
\]
Equating this expression with \eqref{FourthMomentCompressionBound} and rearranging, we obtain
\[
\frac{n-1}{n-2}\frac{f_3^2}{S^2}
+\tr(\mathring{A}^{\perp})^2+\frac{2}{S}|\pi(A\eta_2)|^2
=\frac{n-2}{n}S.
\]
Discarding the last two nonnegative terms yields (2).
Discarding only the last term and using $\frac{n-1}{n-2}\geqslant1$ yields (3).
Thus
\begin{equation}\label{ThirdMomentCompressionBounds}
f_3^2\leqslant\frac{(n-2)^2}{n(n-1)}S^3,\qquad
\frac{f_3^2}{S^2}+\tr(\mathring{A}^{\perp})^2\leqslant\frac{n-2}{n}S.
\end{equation}

For (5), let $\mu_1,\ldots,\mu_{n-2}$ be the eigenvalues of $\mathring{A}^{\perp}$. Then
\[
|\tr(\mathring{A}^{\perp})^3|
=\left|\sum_\alpha\mu_\alpha^3\right|
\leqslant \sum_\alpha|\mu_\alpha|^3
\leqslant \left(\sum_\alpha\mu_\alpha^2\right)^{\frac{3}{2}}=\left(\tr (\mathring{A}^{\perp})^2\right)^{\frac32},
\]
and \eqref{ThirdMomentCompressionBounds} gives
\begin{equation}\label{CenteredCubeAbsorption}
\begin{aligned}
|f_3\tr(\mathring{A}^{\perp})^3|
&\leqslant S\frac{|f_3|}{S}\sqrt{\tr(\mathring{A}^{\perp})^2}\,\tr(\mathring{A}^{\perp})^2\\
&\leqslant \frac{S}{2}\left(\frac{f_3^2}{S^2}+\tr(\mathring{A}^{\perp})^2\right)\tr(\mathring{A}^{\perp})^2\\
&\leqslant \frac{n-2}{2n}S^2\tr(\mathring{A}^{\perp})^2.
\end{aligned}
\end{equation}
\end{proof}

We now return to the estimate of $\mathcal P$~in \eqref{AlgebraicIntegrandDefinition} for $n\geqslant 4$.

\begin{proposition}\label{CurvatureMomentLowerEstimate}
Let $M^n\to\mathbb S^{n+1}(1)$ be a minimal hypersurface with $n\geqslant 4$. At every point where $R\leqslant n(n-1)$, the following estimate holds:
\[
\mathcal P\geqslant \frac{2(n-2)^2}{n^2}S^3(S-n).
\]
  \end{proposition}
 \begin{proof}
First suppose that $S>0$. By \eqref{CurvatureMomentBound}, the curvature assumption is equivalent to $f_4\geqslant S^2-2S$. Hence
  \begin{equation}\label{CurvaturePartOfPowerSum}
\begin{aligned}
\mathcal P
 &=-\left(S-2-\frac{f_4}{S}\right)f_3^2
  +f_3\left(f_5-\frac{f_3f_4}{S}\right)\geqslant  f_3\left(f_5-\frac{f_3f_4}{S}\right).
\end{aligned}
\end{equation} 

Eliminating $f_5$ by \eqref{CenteredCompressionCube} and substituting the expression of $f_4$ from \eqref{FourthMomentCompressionBound}, we obtain
  \begin{align*}
f_3\left(f_5-\frac{f_3f_4}{S}\right)
 &=\frac{1}{3}Sf_3^2
 -\frac{(n-1)(n-3)}{3(n-2)^2}\frac{f_3^4}{S^2}+\frac{f_3^2}{n-2}\tr(\mathring{A}^{\perp})^2
 -\frac{Sf_3}{3}\tr(\mathring{A}^{\perp})^3.
\end{align*} 
Similarly, \eqref{FourthMomentCompressionBound} gives
  \begin{align*}
\frac{(n-2)^2}{n^2}S^2\Bigl(nf_4-(n-2)S^2\Bigr)
&=\frac{(n-2)^2(4-n)}{2n^2}S^4+\frac{(n-2)(n-3)}{2n}Sf_3^2\\
&\quad-\frac{(n-2)^2}{2n}S^3\tr(\mathring{A}^{\perp})^2.
\end{align*} 
Thus,
\begin{align*}
&\quad f_3\left(f_5-\frac{f_3f_4}{S}\right)
 -\frac{(n-2)^2}{n^2}S^2\bigl(nf_4-(n-2)S^2\bigr)\\
&=\left(\frac{(n-2)^2}{n(n-1)}S^3-f_3^2\right)
 \left(\frac{(n-1)(n-3)f_3^2}{3(n-2)^2S^2}
       +\frac{(n-1)(n-4)S}{2n}\right)\\
&\quad+\left(\frac{f_3^2}{n-2}+\frac{(n-2)^2}{2n}S^3\right)
 \tr(\mathring{A}^{\perp})^2-\frac{Sf_3}{3}\tr(\mathring{A}^{\perp})^3\\
&\geqslant\left(\frac{f_3^2}{n-2}+\frac{(n-2)^2}{2n}S^3\right)
 \tr(\mathring{A}^{\perp})^2-\frac{Sf_3}{3}\tr(\mathring{A}^{\perp})^3\\
&\geqslant\left(\frac{f_3^2}{n-2}
 +\frac{(n-2)(3n-7)}{6n}S^3\right)\tr(\mathring{A}^{\perp})^2
 \geqslant0{,}
\end{align*}
where the first inequality follows from $n\geqslant4$ and
\eqref{ThirdMomentCompressionBounds}, and the second follows from
\eqref{CenteredCubeAbsorption}.

Combining these inequalities with \eqref{CurvaturePartOfPowerSum} and using $f_4\geqslant  S^2-2S$, we conclude that
\begin{align*}
f_3f_5-(S-2)f_3^2
&\geqslant \frac{(n-2)^2}{n^2}S^2\bigg(nf_4-(n-2)S^2\bigg)\\
&\geqslant \frac{(n-2)^2}{n^2}S^2
 \bigg(n(S^2-2S)-(n-2)S^2\bigg)\\
&=\frac{2(n-2)^2}{n^2}S^3(S-n).
\end{align*}

At a point where $S=0$, all principal curvatures are zero. Thus both
sides of \eqref{AlgebraicIntegrandTargetBound} vanish, and the inequality
holds with equality.

\end{proof}

We next estimate $\mathcal Q$~in \eqref{DerivativeIntegrandDefinition} to prove \eqref{DerivativeIntegrandTargetBound} on a minimal hypersurface
$M^n\to\mathbb S^{n+1}(1)$ with $n\geqslant4$.
Fix a point $p$
where $S>0$ and choose a principal orthonormal frame.

By \eqref{DerivativeIntegrandDefinition} and
\eqref{CompletedGradientDefinition}, using \eqref{CurvatureMomentDerivatives} gives
\[
\begin{aligned}
\mathcal Q
&=3|T|^2-\frac3{16}|A\nabla S|^2
\leqslant3|T|^2=3\sum_k\left(
\sum_i\left(\lambda_i^2-\frac{\lambda_k\lambda_i}{2}\right)
h_{iik}\right)^2.
\end{aligned}
\]
For each fixed $k$, define a quadratic form on $T_pM$ by
\[
\mathcal Q_k(x)=\left(
\sum_i\left(\lambda_i^2-\frac{\lambda_k\lambda_i}{2}\right)x_i
\right)^2,
\]
where $x=\sum\limits_i x_i e_i$. Then
\begin{equation}\label{DerivativeQuadraticFormDecomposition}
\mathcal Q\leqslant3\sum_k\mathcal Q_k\left(\sum_i h_{iik}e_i\right).
\end{equation}

On the other hand, the Codazzi equation makes $h_{ijk}$ fully symmetric. Separating terms according to the multiplicities of their indices gives
\begin{equation}\label{SymmetricDerivativeMultiplicity}
|\nabla h|^2
=\sum_k\left(3\sum_i h_{iik}^2-2h_{kkk}^2\right)
 +6\sum_{i<j<k}h_{ijk}^2.
\end{equation}
Taking $r=2$ in \eqref{CurvatureMomentDerivatives} and using
\eqref{SymmetricDerivativeMultiplicity}, we also obtain
\begin{equation}\label{DerivativeBlockEnergy}
\begin{aligned}
S^2|\nabla h|^2+S|\nabla S|^2
&\geqslant S^2\sum_k\left(3\sum_i h_{iik}^2-2h_{kkk}^2\right)
 +4S\sum_k\left(\sum_i\lambda_i h_{iik}\right)^2.
\end{aligned}
\end{equation}
Define a quadratic form on $T_pM$ by
\[
\mathcal B_k(x)=S^2\left(3\sum_i x_i^2-2x_k^2\right)
 +4S\left(\sum_i\lambda_ix_i\right)^2.
\]
The same identities give
\[
S^2|\nabla h|^2+S|\nabla S|^2
\geqslant\sum_k\mathcal B_k\left(\sum_i h_{iik}e_i\right).
\]

For each $k$, minimality gives $\sum\limits_i h_{iik}=0$, so the vector $\sum\limits_i h_{iik}e_i$ lies in the subspace
\[
W=\left\{\sum_i x_i e_i\in T_pM\mid\sum_i x_i=0\right\}
=\{x\in T_pM\mid\langle x,\eta_1\rangle=0\}.
\]
Here $\eta_1$ is defined in \eqref{PrincipalFrameVectors}.
To establish \eqref{DerivativeIntegrandTargetBound}, we therefore compare the quadratic forms $\mathcal Q_k$ and $\mathcal B_k$ on $W$, for each $k$.

By \eqref{PrincipalFrameVectors},
$\eta_2=\sum\limits_i\lambda_i e_i$ and
$A\eta_2=\sum\limits_i\lambda_i^2e_i$. Consequently,
\[
\mathcal Q_k(x)
=\left\langle A\eta_2-\frac{\lambda_k}{2}\eta_2,x\right\rangle^2
\qquad(x\in T_pM).
\]
Let $X_k$ be the orthogonal projection of
$A\eta_2-\frac{\lambda_k}{2}\eta_2$ onto $W$. Since
$\langle A\eta_2,\eta_1\rangle=S$ and
$\langle\eta_2,\eta_1\rangle=0$, we have
\[
X_k=A\eta_2-\frac{\lambda_k}{2}\eta_2-\frac{S}{n}\eta_1.
\]
For $x\in W$, the component along $\eta_1$ has zero inner product
with $x$, so
\[
\mathcal Q_k(x)=\langle X_k,x\rangle^2.
\]

Similarly,
\[
\begin{aligned}
\mathcal B_k(x)
&=S^2\left(3\sum_i x_i^2-2x_k^2\right)
 +4S\left(\sum_i\lambda_ix_i\right)^2\\
&=3S^2|x|^2-2S^2\langle e_k,x\rangle^2
 +4S\langle\eta_2,x\rangle^2\\
&=\left\langle 3S^2x-2S^2\langle e_k,x\rangle e_k
 +4S\langle\eta_2,x\rangle\eta_2,\,x\right\rangle.
\end{aligned}
\]
For $x\in W$, define $B_kx$ to be the orthogonal projection onto $W$ of
\[
3S^2x-2S^2\langle e_k,x\rangle e_k+4S\langle\eta_2,x\rangle\eta_2{.}
\]
Then
\[
B_kx=3S^2x-2S^2\langle\bar e_k,x\rangle\bar e_k
 +4S\langle\eta_2,x\rangle\eta_2,
\]
where $\bar e_k=e_k-\frac1n\eta_1$ is the orthogonal projection of $e_k$ onto $W$. Since $\eta_2\in W$ and $\langle\bar e_k,x\rangle=\langle e_k,x\rangle$ for $x\in W$, the operator $B_k:W\to W$ is self-adjoint and satisfies
\[
\mathcal B_k(x)=\langle B_kx,x\rangle\qquad(x\in W).
\]

We now show that $B_k$ is positive definite on $W$.
By Cauchy-Schwarz,
\[
\langle\bar e_k,x\rangle^2
\leqslant|\bar e_k|^2|x|^2=\frac{n-1}{n}|x|^2.
\]
The term $4S\langle\eta_2,x\rangle^2$ is nonnegative. Hence, for every nonzero $x\in W$,
\[
\langle B_kx,x\rangle\geqslant
S^2\left(3-\frac{2(n-1)}{n}\right)|x|^2
=\frac{n+2}{n}S^2|x|^2>0.
\]
Thus $B_k$ is positive definite on $W$.
Let $B_k^{\frac{1}{2}}$ be its positive
square root and $B_k^{-\frac{1}{2}}$ the inverse of this square root.
Applying Cauchy--Schwarz to $B_k^{-\frac12}X_k$ and $B_k^{\frac12}x$, we obtain, for every $x\in W$,
\begin{equation}\label{DerivativeMetricCauchySchwarz}
\begin{aligned}
\mathcal Q_k(x)=\langle X_k,x\rangle^2
&=\langle B_k^{-\frac12}X_k,B_k^{\frac12}x\rangle^2\\
&\leqslant\langle B_k^{-1}X_k,X_k\rangle\langle B_kx,x\rangle
=\langle B_k^{-1}X_k,X_k\rangle\,\mathcal B_k(x).
\end{aligned}
\end{equation}

Consequently, the derivative estimate~\eqref{DerivativeIntegrandTargetBound} reduces to a bound for
$\langle B_k^{-1}X_k,X_k\rangle$.

\begin{lemma}\label{DerivativeDualNormBound}
Let $M^n\to\mathbb S^{n+1}(1)$ be a minimal hypersurface with
$n\geqslant4$. At a point where $S>0$, fix a principal orthonormal frame and $k\in\{1,\ldots,n\}$, and let $X_k$ and $B_k$ be defined as above.
Then
\[
\langle B_k^{-1}X_k,X_k\rangle\leqslant\frac{5n-8}{14(n+2)}.
\]
\end{lemma}
\begin{proof}
Set $Y_k=B_k^{-1}X_k$. Taking inner products of $B_kY_k=X_k$
with $\bar e_k$ and $\eta_2$ gives
\[
S^2
\begin{pmatrix}
\frac{n+2}{n}&\frac{4\lambda_k}{S}\\
-2\lambda_k&7
\end{pmatrix}
\begin{pmatrix}\langle\bar e_k,Y_k\rangle\\\langle\eta_2,Y_k\rangle\end{pmatrix}
=
\begin{pmatrix}\frac{\lambda_k^2}{2}-\frac{S}{n}\\f_3-\frac{S\lambda_k}{2}\end{pmatrix}{,}
\]
where we used
\[
\langle\bar e_k,\eta_2\rangle=\lambda_k,\qquad
\langle\bar e_k,X_k\rangle=\frac{\lambda_k^2}{2}-\frac{S}{n},\qquad
\langle\eta_2,X_k\rangle=f_3-\frac{S\lambda_k}{2}.
\]
Solving this system yields
\[
\begin{aligned}
\langle\bar e_k,Y_k\rangle
&=\frac{n}{S(7(n+2)S+8n\lambda_k^2)}
 \left(\frac{11}{2}\lambda_k^2-\frac{7S}{n}
                  -\frac{4\lambda_k f_3}{S}\right),\\
\langle\eta_2,Y_k\rangle
&=\frac{n}{S(7(n+2)S+8n\lambda_k^2)}
 \left(\lambda_k^3+\frac{n+2}{n}f_3
                  -\frac{n+6}{2n}S\lambda_k\right).
\end{aligned}
\]
Taking the inner product of $B_kY_k=X_k$ with $X_k$, we also have
\[
\begin{aligned}
3S^2\langle Y_k,X_k\rangle
 &=|X_k|^2+2S^2\langle\bar e_k,Y_k\rangle
                    \left(\frac{\lambda_k^2}{2}-\frac{S}{n}\right)
                    -4S\langle\eta_2,Y_k\rangle
                    \left(f_3-\frac{S\lambda_k}{2}\right).
\end{aligned}
\] 
Substituting the two inner products gives
\begin{equation}\label{DerivativeDualNormNumerator}
\begin{aligned}
&3S^2(7(n+2)S+8n\lambda_k^2)\langle Y_k,X_k\rangle\\
=&7(n+2)Sf_4+8n\lambda_k^2f_4-4(n+2)f_3^2
 +(10-3n)S\lambda_kf_3\\
&-16n\lambda_k^3f_3
 +\left(\frac{3n}{4}-\frac{57}{2}\right)S^2\lambda_k^2
 +\frac{19n}{2}S\lambda_k^4-7S^3.
\end{aligned}
\end{equation}
To control the term $-16n\lambda_k^3f_3$, apply Cauchy--Schwarz
to the projections of $e_k$ and $A\eta_2$ onto $V^\perp$ at the
same point. By the definition of $\pi$,
\[
\pi(e_k)=e_k-\frac{1}{n}\eta_1-\frac{\lambda_k}{S}\eta_2,\qquad
\pi(A\eta_2)=A\eta_2-\frac{S}{n}\eta_1-\frac{f_3}{S}\eta_2.
\]
Therefore,
\begin{gather*}
|\pi(e_k)|^2=\frac{n-1}{n}-\frac{\lambda_k^2}{S},\qquad
|\pi(A\eta_2)|^2=f_4-\frac{S^2}{n}-\frac{f_3^2}{S},\\
\langle\pi(e_k),\pi(A\eta_2)\rangle
=\lambda_k^2-\frac{S}{n}-\frac{\lambda_kf_3}{S}.
\end{gather*}
The Cauchy-Schwarz inequality gives
\[
\left(\frac{n-1}{n}-\frac{\lambda_k^2}{S}\right)
\left(f_4-\frac{S^2}{n}-\frac{f_3^2}{S}\right)
\geqslant \left(\lambda_k^2-\frac{S}{n}-\frac{\lambda_kf_3}{S}\right)^2.
\]
Rearranging this inequality gives
\begin{equation}\label{ProjectedCauchySchwarzPolynomial}
\begin{aligned}
\frac{2}{S}f_3\lambda_k^3
 &\geqslant \frac{S^2}{n}-\frac{n-1}{n}f_4
  +\frac{n-1}{nS}f_3^2+\frac{2}{n}\lambda_kf_3-\frac{3S}{n}\lambda_k^2
  +\frac{\lambda_k^2f_4}{S}+\lambda_k^4.
\end{aligned}
\end{equation}
Multiplying \eqref{ProjectedCauchySchwarzPolynomial} by $-8nS$
and applying it to
$-16n\lambda_k^3f_3=-8nS\left(\frac{2}{S}f_3\lambda_k^3\right)$
in \eqref{DerivativeDualNormNumerator} gives
\[
\begin{aligned}
&3S^2(7(n+2)S+8n\lambda_k^2)\langle Y_k,X_k\rangle\\
&\quad\leqslant 3(5n+2)Sf_4-12nf_3^2-3(n+2)S\lambda_kf_3+\frac{3(n-6)}{4}S^2\lambda_k^2
 +\frac{3n}{2}S\lambda_k^4-15S^3.
\end{aligned}
\]
Applying \eqref{FourthMomentCompressionBound} to the $f_4$ term, we obtain
\[
\begin{aligned}
&3S^2(7(n+2)S+8n\lambda_k^2)\langle Y_k,X_k\rangle\\
&\quad\leqslant \frac{3(5n-8)}{2}S^3
 -\frac{9(n^2-n+2)}{2(n-2)}f_3^2-3(n+2)S\lambda_kf_3+\frac{3(n-6)}{4}S^2\lambda_k^2
 +\frac{3n}{2}S\lambda_k^4.
\end{aligned}
\]
For $n\geqslant 4$, the coefficient of $f_3^2$ satisfies
\[
\frac{9(n^2-n+2)}{2(n-2)}-3(n+2)
=\frac{3(n^2-3n+14)}{2(n-2)}>0.
\]
Consequently, completing the square gives
\[
\begin{aligned}
\frac{9(n^2-n+2)}{2(n-2)}f_3^2+3(n+2)S\lambda_kf_3&\geqslant 3(n+2)\bigl(f_3^2+S\lambda_kf_3\bigr)\\
&=3(n+2)\left(
 \left(f_3+\frac{S\lambda_k}{2}\right)^2
 -\frac{S^2\lambda_k^2}{4}\right)\\
&\geqslant -\frac{3(n+2)}{4}S^2\lambda_k^2.
\end{aligned}
\]
The trace-free condition also gives
\[
\lambda_k^2=\left(\sum_{i\ne k}\lambda_i\right)^2
\leqslant (n-1)(S-\lambda_k^2),\qquad
\lambda_k^4\leqslant \frac{n-1}{n}S\lambda_k^2.
\]
Combining these estimates, we obtain
\[
3S^2(7(n+2)S+8n\lambda_k^2)\langle Y_k,X_k\rangle
\leqslant \frac{3(5n-8)}{2}S^3+\frac{3(2n-3)}{2}S^2\lambda_k^2.
\]

It follows that
\[
\begin{aligned}
\langle B_k^{-1}X_k,X_k\rangle
 &=\langle Y_k,X_k\rangle\\
 &\leqslant \frac{\frac{5n-8}{2}S+\frac{2n-3}{2}\lambda_k^2}{7(n+2)S+8n\lambda_k^2}\\
 &=\frac{5n-8}{14(n+2)}
 -\frac{(26n^2-71n+42)\lambda_k^2}
 {14(n+2)(7(n+2)S+8n\lambda_k^2)}\\
 &\leqslant \frac{5n-8}{14(n+2)}.
\end{aligned}
\]
\end{proof}

\begin{proposition}\label{GradientUpperEstimate}
Let $M^n\to\mathbb S^{n+1}(1)$ be a minimal hypersurface with $n\geqslant 4$. Then
\[
\mathcal Q\leqslant \frac{3(5n-8)}{14(n+2)}
\bigl(S^2|\nabla h|^2+S|\nabla S|^2\bigr).
\]
\end{proposition}

\begin{proof}
At a point where $S>0$, apply Lemma~\ref{DerivativeDualNormBound}
for each $k$ in \eqref{DerivativeMetricCauchySchwarz}.
Summing over $k$ in \eqref{DerivativeQuadraticFormDecomposition} and using \eqref{DerivativeBlockEnergy}
yields
 \[
\begin{aligned}
\mathcal Q
&\leqslant \frac{3(5n-8)}{14(n+2)}
 \left(S^2|\nabla h|^2+S|\nabla S|^2
       -6S^2\sum_{i<j<k}h_{ijk}^2\right)\\
&\leqslant \frac{3(5n-8)}{14(n+2)}
       \bigl(S^2|\nabla h|^2+S|\nabla S|^2\bigr).
\end{aligned}
\]
At a point where $S=0$, all $h_{ij}$ vanish. Equation~\eqref{CurvatureMomentDerivatives} then gives $\nabla S=\nabla f_3=0$, so both sides of the desired estimate are zero.
\end{proof}

\begin{theorem}\label{HigherDimensionalRigidity}
Let $M^n\to\mathbb S^{n+1}(1)$ be a closed connected minimal
Tangent-Yang-Mills hypersurface, where $n\geqslant 4$.
If $R\leqslant  n(n-1)$ pointwise on $M$, then $M$ is totally geodesic.
\end{theorem}

\begin{proof}
Propositions~\ref{CurvatureMomentLowerEstimate}
and~\ref{GradientUpperEstimate} give the pointwise estimates
\eqref{AlgebraicIntegrandTargetBound} and
\eqref{DerivativeIntegrandTargetBound}, respectively. These imply
\eqref{HigherDimensionalIntegralComparison}, which together with
\eqref{HigherDimensionalCoefficientGap} gives
 \[
\int_M\bigl(S^2|\nabla h|^2+S|\nabla S|^2\bigr)\,\mathrm{dvol}=0,
\]
and continuity and nonnegativity of the integrand give
\[
S^2|\nabla h|^2+S|\nabla S|^2\equiv0.
\]

Where $S>0$, this identity implies $\nabla h=0$ and hence $\nabla S=0$. Where $S=0$, all $h_{ij}$ vanish, so
$S_k=2\sum\limits_{i,j}h_{ij}h_{ijk}=0$ as well. Therefore $S$ is constant on the connected hypersurface $M$.
Together with the preceding vanishing identity, this gives
$\nabla h\equiv0$. Simons' identity~\eqref{SimonsIdentity} then yields $S(n-S)=0$, so either $S\equiv0$ or $S\equiv n$.

Suppose that $S=n$. The equality case of the classical pinching theorem
\cite{ChernDoCarmoKobayashi} implies that the image of $M$ is a minimal Clifford hypersurface. By the classification of minimal isoparametric Tangent-Yang-Mills hypersurfaces in \cite{GeXiaoZhangHypersurfaces}, the two sphere factors must have the same dimension. Thus $n$ is even
and the image is
\[
\mathbb S^{\frac{n}{2}}\left(\frac{1}{\sqrt{2}}\right)
\times\mathbb S^{\frac{n}{2}}\left(\frac{1}{\sqrt{2}}\right).
\]
Its principal curvatures are $1$ and $-1$, each with multiplicity
$\frac n2$. Consequently,
\[
S=f_4=n,\qquad
R=n(n-1)+n^2-2n-n=2n(n-2)>n(n-1),
\]
where the strict inequality uses $n\geqslant4$. This contradicts the curvature assumption. We conclude that $S\equiv0$ and $h\equiv0$, so $M$ is totally geodesic.
\end{proof}

 \section{The case \texorpdfstring{$n=3$}{n = 3}}\label{DimensionThree}

Let $M^3\to\mathbb S^4(1)$ be a minimal hypersurface, and let all indices
range over $\{1,2,3\}$. 

By minimality, $f_1=0$ and $f_2=S$. Newton's identities give
\[
\sum_{i<j}\lambda_i\lambda_j=-\frac{S}{2},\qquad
\lambda_1\lambda_2\lambda_3=\frac{f_3}{3},
\]
and
\[
f_4=\frac{S}{2}f_2+\frac{f_3}{3}f_1=\frac{S^2}{2},\qquad
f_5=\frac{S}{2}f_3+\frac{f_3}{3}f_2=\frac{5}{6}Sf_3.
\]
The Gauss equation and $f_4=\frac{S^2}{2}$ give
\[
R=6-2S+\frac{S^2}{2}.
\]
Consequently,
\begin{equation}\label{ThreeDimensionalCurvatureRange}
R\leqslant 6\quad\Longleftrightarrow\quad S(S-4)\leqslant 0
\quad\Longleftrightarrow\quad 0\leqslant  S\leqslant 4.
\end{equation}
Suppose that $M$ is closed and satisfies the Tangent-Yang-Mills
equation~\eqref{TangentYangMillsEquation}. Recall from
\eqref{AlgebraicIntegrandDefinition}, \eqref{DerivativeIntegrandDefinition},
and \eqref{CubicMomentIntegralSplit} that
\[
\begin{aligned}
&\mathcal P=f_3f_5-(S-2)f_3^2,\\
&\mathcal Q=\frac13|\nabla f_3|^2
-\frac12\langle A\nabla S,\nabla f_3\rangle,\\
&\int_M\mathcal P\,\mathrm{dvol}=\int_M\mathcal Q\,\mathrm{dvol}.
\end{aligned}
\]
Substituting $f_5=\frac56Sf_3$ gives
\[
\mathcal P=\left(2-\frac S6\right)f_3^2
\]
and hence
\begin{equation}\label{CubicMomentTangentYangMillsIntegral}
\int_M\left(2-\frac S6\right)f_3^2\,\mathrm{dvol}
=\int_M\left(\frac13|\nabla f_3|^2
-\frac12\langle A\nabla S,\nabla f_3\rangle\right)\,\mathrm{dvol}.
\end{equation}

The comparison used for $n\geqslant4$ requires a pointwise lower bound for $\mathcal P$ by a positive multiple of $S^3(S-n)$. In dimension three, the curvature hypothesis gives $0\leqslant S\leqslant4$ and does not provide the corresponding lower bound. Indeed, for every
$3<S\leqslant4$, the trace-free principal curvatures
\[
\left(\sqrt{\frac S2},0,-\sqrt{\frac S2}\right)
\]
satisfy the algebraic curvature bound and have $f_3=0$. Thus
$\mathcal P=0$ while $S^3(S-3)>0$, so no positive constant can give an inequality of the form \eqref{AlgebraicIntegrandTargetBound} for all permitted principal curvatures. This is the obstruction to applying the preceding
comparison directly when $n=3$.

To refine \eqref{CoarseCompletedGradientEstimate}, we retain the principal-curvature weights in $T$.
Applying Cauchy-Schwarz with one factor of $\lambda_i$ separated from the coefficient of $h_{iik}$ gives
\[
\begin{aligned}
|T|^2
&=\sum_k\left(\sum_i\lambda_i
 \left(\lambda_i-\frac{\lambda_k}{2}\right)h_{iik}\right)^2\\
&\leqslant S\sum_{i,k}
 \left(\lambda_i-\frac{\lambda_k}{2}\right)^2h_{iik}^2\\
&\leqslant S\sum_{i,j,k}
 \left(\lambda_i-\frac{\lambda_k}{2}\right)^2h_{ijk}^2\\
&=\frac S4\sum_{i,j,k}
 (5\lambda_i^2-4\lambda_i\lambda_j)h_{ijk}^2.
\end{aligned}
\]

We recall the Peng-Terng invariants~\cite{PengTerng}:
\begin{equation}\label{ThreeDimensionalDerivativeMoments}
\mathcal A=\sum_{i,j,k}\lambda_i^2h_{ijk}^2,
\qquad
\mathcal B=\sum_{i,j,k}\lambda_i\lambda_jh_{ijk}^2.
\end{equation}
The estimate therefore reads
{
\begin{equation*}
|T|^2\leqslant\frac S4(5\mathcal A-4\mathcal B).
\end{equation*}
}

The occurrence of $\mathcal A$ and $\mathcal B$ suggests using a differential identity involving these same derivative terms.
{
The calculation leading to Theorem~4 of Peng and Terng
\cite{PengTerng} gives the following identity in divergence form:}
\begin{equation}\label{CurvatureDerivativeDivergenceIdentity}
\frac12\Div(A^2\nabla S)-\frac13\Div(A\nabla f_3)
=\mathcal A-2\mathcal B+\frac14|\nabla S|^2-Sf_4+f_3^2+S^2.
\end{equation}

As in the proof for $n\geqslant4$, we will combine the
Tangent-Yang-Mills integral identity~\eqref{CubicMomentTangentYangMillsIntegral}
with weighted forms of Simons' identity~\eqref{SimonsIdentity} and
the divergence identity~\eqref{CurvatureDerivativeDivergenceIdentity},
and estimate the resulting zeroth-order terms and first-order derivative terms separately. For smooth functions $f$ and $g$ on an interval containing the range of $S$, \eqref{SimonsIdentity} and
\eqref{CurvatureDerivativeDivergenceIdentity} give the weighted integral identities
\begin{equation}\label{ThreeDimensionalWeightedSimonsIntegral}
0=\int_M\left(
f(S)|\nabla h|^2+\frac{f'(S)}2|\nabla S|^2
+f(S)S(3-S)\right)\,\mathrm{dvol},
\end{equation}
and
\begin{equation}\label{WeightedCurvatureDerivativeIntegral}
\begin{aligned}
0=\int_M\Bigg(&g(S)(\mathcal A-2\mathcal B)
+\frac{g(S)}4|\nabla S|^2
+\frac{g'(S)}2|A\nabla S|^2\\
&-\frac{g'(S)}3\langle A\nabla S,\nabla f_3\rangle
+g(S)\left(f_3^2+S^2-\frac{S^3}{2}\right)\Bigg)\,\mathrm{dvol}.
\end{aligned}
\end{equation}

Summing \eqref{CubicMomentTangentYangMillsIntegral},
\eqref{ThreeDimensionalWeightedSimonsIntegral}, and
\eqref{WeightedCurvatureDerivativeIntegral} yields
\begin{equation}\label{ThreeDimensionalWeightedCombination}
\begin{aligned}
0={}&\int_M\left(\left(2-\frac S6+g(S)\right)f_3^2
+g(S)\left(S^2-\frac{S^3}{2}\right)+f(S)S(3-S)\right)\,\mathrm{dvol}\\
&+\int_M\Bigg(f(S)|\nabla h|^2+g(S)(\mathcal A-2\mathcal B)
+\left(\frac{f'(S)}2+\frac{g(S)}4\right)|\nabla S|^2\\
&\qquad+\frac{g'(S)}2|A\nabla S|^2
+\left(\frac12-\frac{g'(S)}3\right)
 \langle A\nabla S,\nabla f_3\rangle
-\frac13|\nabla f_3|^2\Bigg)\,\mathrm{dvol}.
\end{aligned}
\end{equation}

We choose $f$ and $g$ to be polynomials of degree at most one so that the zeroth-order expression in \eqref{ThreeDimensionalWeightedCombination}
vanishes identically at $S=4$.

\begin{lemma}\label{ThreeDimensionalIntegralDecomposition}
Let $M^3\to\mathbb S^4(1)$ be a closed minimal Tangent-Yang-Mills hypersurface. Define
\begin{align}
&\mathcal L=(4-S)\left(-\frac{f_3^2}{3}
+\frac{S^3}{4}-\frac{S^2}{6}+S\right),
\label{ThreeDimensionalZerothOrderDefinition}\\
&\begin{aligned}
\mathcal M={}&\left(S+\frac43\right)|\nabla h|^2
+\left(\frac S2-\frac{10}{3}\right)(\mathcal A-2\mathcal B)\\
&+\left(\frac S8-\frac13\right)|\nabla S|^2
+\frac13|A\nabla S|^2
-\frac13\left|\nabla f_3-\frac12A\nabla S\right|^2.
\end{aligned}\label{ThreeDimensionalDerivativeDefinition}
\end{align}
Then
\[
0=\int_M(\mathcal L+\mathcal M)\,\mathrm{dvol}.
\]
If $0\leqslant S\leqslant4$, then
\[
\mathcal L\geqslant0,
\]
and $\mathcal L=0$ at a point if and only if $S=0$ or $S=4$.
\end{lemma}
\begin{proof}
To obtain the factor $4-S$, we evaluate the zeroth-order expression
in \eqref{ThreeDimensionalWeightedCombination} at $S=4$:
\[
\left(\frac43+g(4)\right)f_3^2-16g(4)-4f(4).
\]
Requiring this to vanish identically as a polynomial in $f_3^2$ gives
\[
g(4)=-\frac43,\qquad f(4)=\frac{16}{3}.
\]
We therefore write, for constants $c_1,c_2\in\mathbb R$,
\[
f(S)=\frac{16}{3}+c_1(S-4),\qquad
g(S)=-\frac43+c_2(S-4).
\]
The zeroth-order integrand in
\eqref{ThreeDimensionalWeightedCombination} becomes
\[
(4-S)\left(
\left(\frac16-c_2\right)f_3^2
+S\left(4-3c_1+\left(c_1-c_2-\frac23\right)S
+\frac{c_2}{2}S^2\right)\right).
\]
By Lemma~\ref{CompressionMomentInequalities}(2), applied with $n=3$, we have $0\leqslant f_3^2\leqslant\frac{S^3}{6}$.
We require the preceding expression to be nonnegative throughout $0\leqslant S\leqslant4$, for every $f_3^2$ in this interval. For each fixed $S$, its minimum is attained at one of the two endpoints. This requirement is therefore equivalent to
\begin{equation}\label{ThreeDimensionalWeightEndpointConditions}
\left\{
\begin{aligned}
&4-3c_1+\left(c_1-c_2-\frac23\right)S
+\frac{c_2}{2}S^2\geqslant0,\\
&4-3c_1+\left(c_1-c_2-\frac23\right)S
+\left(\frac{c_2}{3}+\frac1{36}\right)S^2\geqslant0.
\end{aligned}
\right.
\end{equation}
At $S=0$, these polynomial conditions require
\[
c_1\leqslant\frac43.
\]
One convenient choice is
\[c_1=1,\qquad c_2=\frac12.\]
which satisfy \eqref{ThreeDimensionalWeightEndpointConditions}.
Now the functions $f(S)$ and $g(S)$ in
\eqref{ThreeDimensionalWeightedCombination} are chosen as
\[
f(S)=S+\frac43,\qquad g(S)=\frac S2-\frac{10}{3}.
\]
For this choice, using the definition of $\mathcal L$ in
\eqref{ThreeDimensionalZerothOrderDefinition} and the bound
$f_3^2\leqslant\frac{S^3}{6}$ gives
\[
\mathcal L\geqslant\frac{(4-S)S}{36}(7S^2-6S+36)\geqslant0.
\]
Equality holds precisely when $S=0$ or $S=4$. Substituting $f(S)=S+\frac43$ and $g(S)=\frac S2-\frac{10}{3}$ into
\eqref{ThreeDimensionalWeightedCombination} and completing the square
gives the integral identity.

\end{proof}

For $i\in\{1,2,3\}$, define $i'$ and $i''$ by

\[
\begin{array}{c|ccc}
i&1&2&3\\\hline
i'&2&3&1\\
i''&3&1&2
\end{array}
\]
At a fixed point, choose an orthonormal principal frame. By
\eqref{CurvatureMomentDerivatives},
\begin{align*}
|\nabla S|^2
&=4\sum_k\left(\sum_j\lambda_jh_{jjk}\right)^2,\\
|A\nabla S|^2
&=4\sum_k\left(\sum_j\lambda_k\lambda_jh_{jjk}\right)^2,\\
\left|\nabla f_3-\frac12A\nabla S\right|^2
&=\sum_k\left(\sum_j(3\lambda_j^2-\lambda_k\lambda_j)h_{jjk}\right)^2.
\end{align*}
The full symmetry of $h_{ijk}$ also gives
\[
|\nabla h|^2
=\sum_i\left(h_{iii}^2+3h_{i'i'i}^2+3h_{i''i''i}^2\right)+6h_{123}^2.
\]
Expanding the sums in \eqref{ThreeDimensionalDerivativeMoments} gives
\begin{align*}
\mathcal A={}&\sum_i\Bigl(\lambda_i^2h_{iii}^2
+(2\lambda_{i'}^2+\lambda_i^2)h_{i'i'i}^2
+(2\lambda_{i''}^2+\lambda_i^2)h_{i''i''i}^2\Bigr)
+2Sh_{123}^2,\\
\mathcal B={}&\sum_i\Bigl(\lambda_i^2h_{iii}^2
+(\lambda_{i'}^2+2\lambda_i\lambda_{i'})h_{i'i'i}^2
+(\lambda_{i''}^2+2\lambda_i\lambda_{i''})h_{i''i''i}^2\Bigr)
-Sh_{123}^2.
\end{align*}
Thus the terms with three distinct indices
in $\mathcal M$ of \eqref{ThreeDimensionalDerivativeDefinition} contribute
\[
\left(6\left(S+\frac43\right)
+4S\left(\frac S2-\frac{10}{3}\right)\right)h_{123}^2
=\frac23(3S^2-11S+12)h_{123}^2.
\]
For each $i$, define
\begin{align*}
\mathcal M_i={}&\left(S+\frac43\right)
\left(h_{iii}^2+3h_{i'i'i}^2+3h_{i''i''i}^2\right)\\
&+\left(\frac S2-\frac{10}{3}\right)
\Bigl(-\lambda_i^2h_{iii}^2
+(\lambda_i^2-4\lambda_i\lambda_{i'})h_{i'i'i}^2
+(\lambda_i^2-4\lambda_i\lambda_{i''})h_{i''i''i}^2\Bigr)\\
&+\left(\frac S2-\frac43+\frac43\lambda_i^2\right)
\left(\sum_j\lambda_jh_{jji}\right)^2
-\frac13\left(\sum_j(3\lambda_j^2-\lambda_i\lambda_j)h_{jji}\right)^2.
\end{align*}
The preceding identities give an alternative expression for $\mathcal M$ in
\eqref{ThreeDimensionalDerivativeDefinition}:
\begin{equation}\label{ThreeDimensionalDerivativeBlocks}
\mathcal M=\sum_{i=1}^3\mathcal M_i
+\frac23(3S^2-11S+12)h_{123}^2.
\end{equation}

For each fixed $i\in\{1,2,3\}$, set
\[
\widetilde\lambda_i=\lambda_{i'}-\lambda_{i''},\qquad
x=h_{iii},\qquad y=h_{i'i'i}-h_{i''i''i}.
\]
The trace-free conditions give
\[
\begin{gathered}
\lambda_{i'}=\frac{\widetilde\lambda_i-\lambda_i}{2},\qquad
\lambda_{i''}=-\frac{\lambda_i+\widetilde\lambda_i}{2},\qquad
S=\frac{3\lambda_i^2+\widetilde\lambda_i^2}{2},\\
h_{iii}=x,\qquad
h_{i'i'i}=\frac{y-x}{2},\qquad
h_{i''i''i}=-\frac{x+y}{2}.
\end{gathered}
\]
Consequently,
\begin{align*}
h_{iii}^2+3h_{i'i'i}^2+3h_{i''i''i}^2
&=x^2+\frac34(y-x)^2+\frac34(x+y)^2
=\frac52x^2+\frac32y^2,\\
\sum_j\lambda_jh_{jji}
&=\lambda_i x+\frac14(\widetilde\lambda_i-\lambda_i)(y-x)
+\frac14(\lambda_i+\widetilde\lambda_i)(x+y)=\frac32\lambda_i x+\frac12\widetilde\lambda_i y.
\end{align*}
For the contribution of $\mathcal A-2\mathcal B$, substitution gives
\begin{align*}
&-\lambda_i^2h_{iii}^2
+(\lambda_i^2-4\lambda_i\lambda_{i'})h_{i'i'i}^2
+(\lambda_i^2-4\lambda_i\lambda_{i''})h_{i''i''i}^2\\
=&-\lambda_i^2x^2
+\frac14(3\lambda_i^2-2\lambda_i\widetilde\lambda_i)(y-x)^2
+\frac14(3\lambda_i^2+2\lambda_i\widetilde\lambda_i)(x+y)^2\\
=&\;\frac{\lambda_i^2}{2}x^2
+2\lambda_i\widetilde\lambda_i xy+\frac{3\lambda_i^2}{2}y^2.
\end{align*}
On the other hand,
\begin{align*}
3\lambda_{i'}^2-\lambda_i\lambda_{i'}
&=\frac34(\widetilde\lambda_i-\lambda_i)^2
-\frac12\lambda_i(\widetilde\lambda_i-\lambda_i)=\frac14(5\lambda_i^2-8\lambda_i\widetilde\lambda_i
+3\widetilde\lambda_i^2),\\
3\lambda_{i''}^2-\lambda_i\lambda_{i''}
&=\frac34(\lambda_i+\widetilde\lambda_i)^2
+\frac12\lambda_i(\lambda_i+\widetilde\lambda_i)=\frac14(5\lambda_i^2+8\lambda_i\widetilde\lambda_i
+3\widetilde\lambda_i^2).
\end{align*}
It follows that
\begin{align*}
\sum_j(3\lambda_j^2-\lambda_i\lambda_j)h_{jji}
&=2\lambda_i^2x
+\frac18(5\lambda_i^2-8\lambda_i\widetilde\lambda_i
+3\widetilde\lambda_i^2)(y-x)-\frac18(5\lambda_i^2+8\lambda_i\widetilde\lambda_i
+3\widetilde\lambda_i^2)(x+y)\\
&=\frac34(\lambda_i^2-\widetilde\lambda_i^2)x
-2\lambda_i\widetilde\lambda_i y.
\end{align*}
Substituting these expansions, we obtain
\begin{equation}\label{ThreeDimensionalBlockExpansion}
\begin{aligned}
\mathcal M_i={}&\left(S+\frac43\right)
\left(\frac52x^2+\frac32y^2\right)+\left(\frac S2-\frac{10}{3}\right)
\left(\frac{\lambda_i^2}{2}x^2
+2\lambda_i\widetilde\lambda_i xy+\frac{3\lambda_i^2}{2}y^2\right)\\
&+\left(\frac S8-\frac13+\frac{\lambda_i^2}{3}\right)
(3\lambda_i x+\widetilde\lambda_i y)^2-\frac13\left(\frac34(\lambda_i^2-\widetilde\lambda_i^2)x
-2\lambda_i\widetilde\lambda_i y\right)^2.
\end{aligned}
\end{equation}
To establish the positivity of these blocks, we first prove an algebraic lemma.

\begin{lemma}\label{ThreeDimensionalQuadraticFormPositivity}
Let $a,b\in\mathbb R$ satisfy $3a^2+b^2\leqslant8$, and define
\begin{align*}
\alpha(a,b)&=234a^4+51a^2b^2-9b^4-44a^2+60b^2+160,\\
\beta(a,b)&=ab(135a^2-3b^2-208),\\
\gamma(a,b)&=54a^4-21a^2b^2+3b^4-132a^2+20b^2+96.
\end{align*}
Then
\[
\alpha(a,b)>0,\qquad
\alpha(a,b)\gamma(a,b)-\beta(a,b)^2>0.
\]
In particular, the quadratic form
\[
\alpha(a,b)x^2+2\beta(a,b)xy+\gamma(a,b)y^2
\]
is positive definite in $(x,y)$.
\end{lemma}
\begin{proof}
In this proof, we write $\alpha,\beta,\gamma$ for
$\alpha(a,b),\beta(a,b),\gamma(a,b)$, respectively. Set
\[
t=\sqrt{4-\frac{3a^2+b^2}{2}}.
\]
Then
\[
3a^2+b^2+2t^2=8.
\]
Using this identity to homogenize the coefficients gives
\begin{align*}
\alpha
=&234a^4+51a^2b^2-9b^4
+\frac12(3a^2+b^2+2t^2)(-11a^2+15b^2)+\frac52(3a^2+b^2+2t^2)^2\\
=&240a^4+83a^2b^2+b^4+19a^2t^2+25b^2t^2+10t^4.
\end{align*}
Similarly,
\begin{align*}
\gamma
=&54a^4-21a^2b^2+3b^4
+\frac12(3a^2+b^2+2t^2)(-33a^2+5b^2)+\frac32(3a^2+b^2+2t^2)^2\\
=&18a^4-21a^2b^2+7b^4-15a^2t^2+11b^2t^2+6t^4,
\end{align*}
and
\[
\beta=ab\bigl(135a^2-3b^2-26(3a^2+b^2+2t^2)\bigr)
=ab(57a^2-29b^2-52t^2).
\]
Since $a,b,t$ cannot vanish simultaneously,
\[
\alpha\geqslant240a^4+b^4+10t^4>0.
\]
For the determinant, direct expansion gives
\begin{equation}\label{ThreeDimensionalDeterminantSquares}
\begin{aligned}
512(\alpha\gamma-\beta^2)
={}&8(524a^4-435a^2b^2-204a^2t^2
            +21b^4+88b^2t^2+58t^4)^2\\
&+32\bigl(ab(68a^2-5b^2-33t^2)\bigr)^2
+4\bigl(at(88a^2+99b^2-206t^2)\bigr)^2\\
&+\bigl(bt(1296a^2-253b^2-194t^2)\bigr)^2\\
&+15232(a^4)^2+20032(a^3b)^2+1528(a^2b^2)^2\\
&+1488(ab^3)^2+56(b^4)^2+11264(a^3t)^2+12160(a^2bt)^2\\
&+1916(ab^2t)^2+1655(b^3t)^2+9344(a^2t^2)^2+80(abt^2)^2\\
&+108(b^2t^2)^2+1136(at^3)^2+13820(bt^3)^2+3808(t^4)^2\\
\geqslant{}&15232a^8+56b^8+3808t^8>0.
\end{aligned}
\end{equation}
The bound \eqref{ThreeDimensionalDeterminantSquares} and $\alpha>0$
imply positive definiteness.
\end{proof}

\Needspace{8\baselineskip}
\begin{theorem}\label{ThreeDimensionalRigidity}
Let $M^3\to\mathbb S^4(1)$ be a closed connected minimal
Tangent-Yang-Mills hypersurface.
If $R\leqslant6$ pointwise on $M$, then $M$ is totally geodesic.
\end{theorem}

\begin{proof}
By \eqref{ThreeDimensionalCurvatureRange}, we have
$0\leqslant S\leqslant4$.
Lemma~\ref{ThreeDimensionalIntegralDecomposition} gives
\[
0=\int_M\mathcal L\,\mathrm{dvol}
+\int_M\mathcal M\,\mathrm{dvol},
\qquad \mathcal L\geqslant0.
\]

We next show that $\mathcal M\geqslant0$. At any fixed point, choose an orthonormal principal frame.
For each $i$, recall that
$x=h_{iii}$ and $y=h_{i'i'i}-h_{i''i''i}$.
Substituting $S=\frac{(3\lambda_i^2+\widetilde\lambda_i^2)}{2}$
into \eqref{ThreeDimensionalBlockExpansion}
and collecting coefficients yields
\[
\begin{aligned}
48\mathcal M_i=
\alpha(\lambda_i,\widetilde\lambda_i)x^2
+2\beta(\lambda_i,\widetilde\lambda_i)xy+\gamma(\lambda_i,\widetilde\lambda_i)y^2.
\end{aligned}
\]
Since $3\lambda_i^2+\widetilde\lambda_i^2=2S\leqslant8$,
Lemma~\ref{ThreeDimensionalQuadraticFormPositivity}
implies that $\mathcal M_i$ is positive definite in $(x,y)$.
Moreover,
\[
3S^2-11S+12
=3\left(S-\frac{11}{6}\right)^2+\frac{23}{12}>0.
\]
The decomposition \eqref{ThreeDimensionalDerivativeBlocks}
therefore gives $\mathcal M\geqslant0$.

Both integrands in the integral identity are continuous and
nonnegative, so
\[
\mathcal L=\mathcal M=0
\]
everywhere on $M$.
The positive definiteness of the blocks $\mathcal M_i$
and the positivity of the coefficient of $h_{123}^2$
in~\eqref{ThreeDimensionalDerivativeBlocks} then imply
\[
h_{123}=0,\qquad h_{iii}=0,\qquad
h_{i'i'i}-h_{i''i''i}=0
\quad(i=1,2,3).
\]
The trace-free condition
$h_{iii}+h_{i'i'i}+h_{i''i''i}=0$
gives $h_{i'i'i}=h_{i''i''i}=0$ for every $i$.
By the full symmetry of $h_{ijk}$, all components of
$\nabla h$ vanish, so $\nabla h\equiv0$.

Hence $S$ is constant, and Simons' identity gives $S=0$ or $S=3$.
On the other hand, the equality case $\mathcal L=0$ in
Lemma~\ref{ThreeDimensionalIntegralDecomposition}
implies $S=0$ or $S=4$.
Therefore $S=0$, and $M$ is totally geodesic in $\mathbb S^4(1)$.
\end{proof}

{\bf AI disclosure.}
During the preparation of this manuscript, we used AI to improve the language and presentation and to assist in establishing the positive definiteness of the quadratic forms in Lemma 4.2. All AI-assisted mathematical derivations and computations were independently checked by the authors, who reviewed and finalized the manuscript. We take full responsibility for its content.

\end{document}